\documentclass[a4paper,11pt,reqno]{article}
\usepackage[T1]{fontenc}
\usepackage{courier}
\usepackage[scaled=0.92]{helvet}

\usepackage{amscd,amsmath,amssymb,amsthm,amsfonts,epsfig,graphics,accents}

\usepackage[english]{babel}
\usepackage{amsthm,amsfonts,amssymb,amsmath}
\usepackage{hyperref}
\hypersetup{
colorlinks=true,
citecolor=red,
linkbordercolor={1 1 1},
citebordercolor={1 1 1},
}
\usepackage{cleveref}
\usepackage{mathrsfs}

\usepackage{xcolor} 

\theoremstyle{plain}
\newtheorem{satz}{Theorem}[section]
\newtheorem{prop}[satz]{Proposition}
\newtheorem{cor}[satz]{Corollary}
\newtheorem{lem}[satz]{Lemma}

\theoremstyle{definition}

\newtheorem{rem}[satz]{Remark}

\newtheorem{hyp}[satz]{Hypothesis}

\newcommand{\mx}{\mbox}
\newcommand{\rw}{\rightarrow}
\newcommand{\de}{\displaystyle}
\newcommand{\ml}{\mathcal}

\newcommand{\pl}{\partial}

\newcommand{\beq}[1]{\begin{equation} \label{#1}}
\newcommand{\eeq}{\end{equation}}
\newcommand{\beqar}[1]{\beq{#1} \begin{array}{rcl}}
\newcommand{\eeqar}{\end{array} \eeq}

\newcommand{\jj}{^{(j)}}

\newcommand{\hrho}{\hat{\rho}}

\newcommand{\drx}[1]{\check{#1}}
\newcommand{\ddrx}[1]{\check{\check{#1}}}

\providecommand{\ep}{\varepsilon}
\providecommand{\ph}{\varphi}
\providecommand{\ka}{\kappa}

\providecommand{\RR}{\mathbb{R}}
\providecommand{\CC}{\mathbb{C}}

\providecommand{\ZZ}{\mathbb{Z}}
\providecommand{\NN}{\mathbb{N}}

\providecommand{\TT}{\mathbb{T}}

\newcommand{\norm}[2]{\left \lVert#1 \right\rVert_{#2}}

\newcommand{\ff}[1]{(\ref{#1})}

\usepackage[utf8]{inputenc}

\hypersetup{pdfauthor={Fortunati, A.},%
            pdftitle={Travelling waves over an arbitrary bathymetry...},%
            pdfsubject={},%
            pdfkeywords={}%
}

\DeclareMathOperator{\sign}{sign}

\makeatletter
\renewcommand*{\@fnsymbol}[1]{\ensuremath{\ifcase#1\or *\or \mathsection \or (b)\or \else \fi}}
\makeatother

\title{\LARGE{\textbf{Travelling waves over an arbitrary bathymetry: a perturbation theory approach.}}}
\date{}
\author{%
Alessandro Fortunati\thanks{E-mail: alessandro.fortunati@bristol.ac.uk} 
\bigskip \\
School of Mathematics, University of Bristol, Bristol BS8 1TW, United Kingdom.
}

\begin{document}

\maketitle

\begin{abstract}
The problem of a travelling wave over an arbitrary quasi-flat bathymetry in a semi infinite channel is studied in the shallow-water formulation. It is shown how the  streamfunction can be cast, in the vicinity of an elliptic equilibrium for the fluid flow, in the form of a nearly-integrable non-autonomous Hamiltonian with \emph{aperiodic} time dependence. The proofs use the tools of perturbation theory in the real-analytic setting. \\
The obtained Hamiltonian provides a natural example in the context of the aperiodically time-dependent Hamiltonian systems studied in Fortunati and Wiggins (2016). Some key properties of the system at hand, such as the stability, can be addressed as a consequence of this theory.
\smallskip\\
{\it Keywords:} Shallow-water equations, Aperiodic time dependence, Nearly-integrable Hamiltonian systems.
\smallskip\\
{\it 2010 MSC}: Primary: 35C07, 70H08. Secondary 35C10, 37J40. 
 
\end{abstract}

\section{Preliminaries}\label{sec:one}
This paper deals with the problem of a propagating wave for the linear quasi-geostrophic shallow-water equation for an inviscid fluid in terms of the \emph{streamfunction} $\psi=\psi(x,y,t)$,
\beq{eq:sw}
\pl_{t}(\Delta \psi-F \psi)+J(\psi,\delta)+\mathcal{F} \pl_x \psi=0 \mx{,}
\eeq
where $J(f_1,f_2):=\pl_x f_1 \pl_y f_2 - \pl_y f_1 \pl_x f_2$, $F, \ml{F} \in \RR$ and $\delta=\delta(x,y)$, see \ff{eq:bath}.\\
For the derivation of equation \ff{eq:sw} and more details about the presentation below, until (\ref{eq:dispersion}),  we refer to the comprehensive exposition of \cite{psky}. \\
The solutions of \ff{eq:sw} will be studied in a \emph{semi-infinite channel}
\[
\ml{C}:=\RR^+ \times \TT \ni (x,y) \mx{,}
\]
here $\RR^+:=(0,+\infty)$, rotating with angular velocity $\Omega$ around the vertical axis $z$ and with a depth $d$, which is supposed to be ``small'' in a suitable scaling, i.e. $d \ll 1$. The bottom of the channel is assumed to be described by a quasi-flat bathymetry of the form
\beq{eq:bath}
\delta(x,y):=\mu \tilde{g}(x,y) - d \mx{,}
\eeq
in the hypothesis $0\ll \mu \ll d $. The assumptions on $\tilde{g}$ will be discussed later.\\
The parameter $F=\ml{O}(1)$ is related to the geometric scaling of the problem, see \cite[formula 3.12.9]{psky} while $\ml{F}=\ml{O}(1)$ represents the contribution of the linear term in the well known $\beta-$\emph{plane} approximation\footnote{More precisely, we have set $\epsilon:=r_0^{-1}\tan \theta_0$ where $r_0$ is the Earth radius and $\theta_0$ is the mean latitude, see \cite[formula 3.17.5]{psky} then we have chosen $\epsilon \ml{F}:=\ep$, where $\ep \equiv \ep_T \ll 1$, see \cite[Sec. 3.12]{psky}, so that $\epsilon$ and $\ep$ have the same magnitude.}.\\
The equation (\ref{eq:sw}) is linear as the term $ J(\Delta \psi - F \psi,\psi)$ has been neglected. This assumption can be thought as ``appropriate'' in the case of motions generated by small amplitude streamfunctions\footnote{In other terms, this can be seen by considering the rescaling $\psi \rw \delta \psi$ then neglecting the $\ml{O}(\delta^2)$ contributions.}. However, we remark that the linear case itself exhibits some peculiar and non-trivial difficulties due to the presence of a non-periodic bathymetry that will be described in detail later on.   \\
Given a streamfunction $\psi$, the velocity field is associated via the \emph{geostrophic equations}
\beq{eq:geo}
\dot{y}=\pl_x \psi, \qquad \dot{x}=-\pl_y \psi \mx{,}
\eeq
which represent the lowest order approximation in the quasi-geostrophic scaling. These exploit the intrinsic Hamiltonian structure of the problem.\\
The impenetrability of the channel boundaries is expressed by the boundary conditions
\beq{eq:boundary}
\pl_x \psi|_{y=0,2 \pi}  \equiv 0, \qquad \forall x \geq 0 \mx{.} 
\eeq
Clearly, if $\mu=0$ (flat bathymetry), equation (\ref{eq:sw}) is simply reduced to the homogeneous equation 
\beq{eq:swzero}
\pl_t (\Delta \psi^{(0)} -F \psi^{(0)})+\mathcal{F} \pl_x \psi^{(0)} = 0\mx{.}
\eeq
A solution for the previous equation can be chosen in the form of \emph{travelling wave} 
\beq{eq:psizero}
\psi^{(0)}(x,y,t)=A \sin (\tilde{m} y) \cos (\kappa x + \sigma t) \mx{,}
\eeq
where $A\in \RR$ is the \emph{amplitude}, while $(\tilde{m},\kappa) \in \NN^2 \setminus \{(0,0)\}$ and $\sigma \in \RR$ are commonly known as \emph{wave numbers} and \emph{wave speed}, respectively. The parameter $A$ will be supposed fixed once and for all.\\
Equation (\ref{eq:psizero}) is satisfied once $\sigma$ has been chosen as follows (\emph{dispersion relation})
\beq{eq:dispersion}
\sigma=\sigma(\ka) := \frac{\kappa \mathcal{F}}{\tilde{m}^2+\kappa^2+F} \mx{.}
\eeq
Note that $\sigma(-\ka)=-\sigma(\ka)$ and $\sigma(\ka)  \neq 0$.\\
Without loss of generality we shall suppose $\mathcal{F}<0$ in such a way $\sigma < 0$ and the wave propagates in the direction of the positive $x-$axis. \medskip \\
It is customary, when dealing with travelling waves problems, see e.g. \cite{kn87}, to consider the following Galileian transformation
\beq{eq:galileian}
\mathcal{G}: \qquad (y,x)=(p,q-\ka^{-1} \sigma t) \mx{.}
\eeq
The equations of motion for the streamlines given by \ff{eq:psizero} via \ff{eq:geo} in these new variables 
\[
\left\{
\begin{array}{rcl}
\de \dot{p} & = & \de - A \ka \sin (\tilde{m} p) \sin (\ka q)\\
\de \dot{q} & = & \de \frac{\sigma}{\ka} -\de A \tilde{m} \cos (\tilde{m} p) \cos ( \ka q)  
\end{array} \mx{,}
\right.
\]
are generated by the Hamiltonian 
\beq{eq:ht}
\tilde{H}_0:=\frac{\sigma}{\ka} p - A \sin (\tilde{m} p) \cos (\ka q) \mx{.} 
\eeq
Under the condition 
\beq{eq:condparam}
-\sigma \leq \kappa \tilde{m} A \mx{,}
\eeq
it is immediate to check that the system above possesses the following equilibria in $[-\pi,\pi]^2$
\beq{eq:equilibria}
(p_e^{\pm},q_e)=\left(\pm \arccos \left( \frac{\sigma}{\ka \tilde{m} A} \right),0 \right), \qquad
(p_h,q_h^{\pm})=\left(0, \pm \frac{1}{\tilde{m}} \arccos \left( \frac{\sigma}{\ka \tilde{m} A} \right) \right) \mx{,}
\eeq
and that $(p_e^{\pm},q_e)$  and $(p_h,q_h^{\pm})$ are elliptic and hyperbolic, respectively. \medskip \\
Given its remarkable geophysical application, the study of the disturbances on a fluid arising from a ``rough'' bottom, is a widely studied problem and it has been investigated for a long time (see e.g. \cite{rhi73}) in many different settings and regimes. The literature on this topic is remarkably vast and a review goes far beyond the purpose of this paper, see e.g. \cite{van03}, \cite{ch07}, \cite{cr10}, \cite{cr12} and references therein. The paper \cite{cr12}, in particular, contains an extensive analysis from the Hamiltonian point of view, in which relevant cases of periodic and $O(1)$ bottoms are studied. The case of a random bottom, already investigated in \cite{ro83}, has been developed in \cite{cr10}.\\  
More closely to the travelling wave problem, a remarkable attention has been devoted to the Rossby waves, either oceanic or atmospheric. The key interpretation of a travelling wave as an integrable and autonomous (Hamiltonian) system for the streamlines, \cite{kw89}, has led to the study of several classes of disturbances to be interpreted as ``perturbations'', see, for instance, \cite{pi91} and \cite{ma98}. In this setting, phenomena such as \emph{chaotic transport} or \emph{chaotic mixing} have been investigated in the ``perturbed'' systems, either numerically or by means of a Melnikov analysis of the stable and unstable manifolds asymptotic to hyperbolic points like $(p_h,q_h^{\pm})$.\\
In general, the above described studies stress the importance of results of stability in the presence of \emph{aperiodically} time-dependent perturbation as ``barriers'' for the Lagrangian transport in fluid dynamics. See \cite{sam06} for an extensive review of the subject.   \medskip  \\
Starting from special (integrable) model given by the streamfunction of the travelling wave \ff{eq:psizero} in a neighbourhood of the point $(p_e^{+},q_e)$, see $H_0(I)$ in \ff{eq:hamiltonian}, our aim is to construct a perturbation arising from a ``slightly'' variable bottom as described by \ff{eq:bath}. For this purpose, the bathymetry will be supposed to be flat for $x \rw +\infty$ but without further assumptions on its dependence on $x$ (e.g. periodicity). It will be shown that, at least in the vicinity of the point $(p_e^{+},q_e)$, the streamlines are described by a nearly-integrable Hamiltonian system with a non-autonomous perturbation in a suitable real-analyticity class. It might be redundant to stress that the construction of the mentioned perturbation itself, which is basically a result of existence for the solution of the ``perturbed problem'' (i.e. the equation  (\ref{eq:sw}) with $\mu>0$), is the difficult part of the whole result. We anticipate that, despite the linearity of the problem, the general (non-periodic) dependence on $x$ of the bathymetry inhibits the standard diagonalization of the (linear) perturbed operator in the Fourier space. This is a commonly used tool when dealing with the linearised operator in the context of superconvergent methods for non-linear PDEs, see e.g. \cite[Section 4.5]{be03}. The above described invertibility of the linear operator has proven, instead, by means of perturbative methods. \\
Once this has been shown, it is easy to realize that, due to the general (non-periodic) dependence on $x$ of $\tilde{g}$, the transformation $\mathcal{G}$ gives rise to an aperiodically time-dependent perturbation. The described class of Hamiltonian systems has been studied in \cite{fw16} and the case at hand provides a natural example in which the mentioned results can be applied.\\

\section{Analytic setting and main statement}\label{sec:two}
Given $\rho \in (0,1/2]$ we introduce the complexification of $\ml{C}$ by defining
\[
\ml{C}_{\rho}:=\ml{S}_{\rho} \times \TT_{\rho}
\]
where
\[
\ml{S}_{\rho}:=\{x \in \CC: \Re x > -\rho, \quad |\Im x|<\rho \}, \qquad 
\TT_{\rho}:=\{y \in \TT: |\Im y|<\gamma  \} \mx{.}
\]
Given a function $v:C_{\rho} \rw \CC$, we define, for all $\rho'<\rho$, the ($x-$dependent) Fourier norm 
\[
\norm{v}{\rho'}:=\sum_{k \in \ZZ} |v_k(x)| e^{|k|\rho'} \mx{,}
\]
where $v_k(x)$ are the Fourier coefficients $v(x,y):=\sum_{z \in \ZZ} v_k(x)e^{i k y}$.
\\ 
Similarly, for any function $\tilde{v}=\tilde{v}(x,y,t): \ml{C}_{\rho} \times \RR \rw \CC$, the norm $\norm{\tilde{v}}{\rho'}$ will be dependent on $(x,t)$.\\
Let us now consider the following set $\mathcal{U}:=G \times \TT \ni (I,\ph)$ where $G \subset \RR$. Given $\gamma \in (0,1/2]$, the set $\ml{U}$ can be complexified by considering
\[
\ml{U}_{\gamma}:=G_{\gamma} \times \TT_{\gamma}, \qquad  G_{\gamma}:=\{I^* \in \CC: |I^*-I|<\gamma, \quad \forall I \in G \} \mx{,}
\]  
and it can be endowed with the following norm  
\[
\norm{w}{\gamma'}^* := \sum_{k \in \ZZ} |w_k| e^{|k| \gamma'}, \qquad \forall \gamma'<\gamma \mx{,}
\]
where $w_k=w_k(I)$ and $\tilde{w}_k=\tilde{w_k}(I,t)$ if $w:\mathcal{U}_{\gamma} \rw \CC$ or $\tilde{w}=\tilde{w}(I,\ph,t):\mathcal{U}_{\gamma} \times \RR \rw \CC$, respectively. Let us finally define $\ml{W}_{\gamma}:=\{\eta \in \CC : |\Im \eta|<\gamma\}$. \\
Throughout the paper we shall consider the following class of topographies
\begin{hyp}\label{hyp} The function $\tilde{g}:\ml{C} \rw \RR$ is holomorphic on $\ml{C}_{2 \rho}$ and satisfies
\beq{eq:hypg}
\norm{\tilde{g}(x,y)}{\rho} \leq \ml{M} e^{-\nu |x|} \mx{,}
\eeq 
for some\footnote{Note that no lower bounds on $\nu$ are imposed, consistently the threshold for the values of $\mu$ allowed in Theorem \ref{thm}, becomes smaller and smaller with $\nu$, see \ff{eq:mustar}.} $\nu>0$. 
\end{hyp}
In the above described setting, the main result states as follows
\begin{satz}\label{thm} Assume hypothesis \ref{hyp} and condition \ff{eq:condparam}. Then for sufficiently small\footnote{See \ff{eq:mustar} for a quantitative estimate.} $\mu$, there exists $0<\tilde{\gamma}<\gamma$ and a set of coordinates $(I,\ph) \in \mathcal{U}_{\tilde{\gamma}}$ for some suitable interval $G$, such that the solutions of (\ref{eq:sw}), in a neighbourhood of $(p_e^+,q_e)$, can be cast in the form 
\beq{eq:hamiltonian}
H_{\mu}(I,\ph,\eta,t):=H_0(I)+ \eta + H_1(I,\ph,t;\mu) \mx{,}
\eeq
where $H_{\mu}$ is a real-analytic function over $\ml{U}_{\tilde{\gamma}}\times \ml{W}_{\tilde{\gamma}}$ and in $\mu$ for all $t \geq 0$. We denote with $\eta$ the variable canonically conjugated to $t$. Furthermore 
\beq{eq:decayintime}
\sup_{I \in G_{\tilde{\gamma}}}\norm{H_1}{\tilde{\gamma}} \leq K e^{-\tilde{\nu} t} \mx{,}
\eeq
for some $K=\mathcal{O}(\mu)>0$ and $\tilde{\nu}:=\tilde{\nu}(\nu)$ with $\tilde{\nu}(0)=0$. 
\end{satz}
Similarly to the aperiodic dependence on $t$ of the perturbation, the time decay \ff{eq:decayintime} is another direct consequence of the transformation $\mathcal{G}$ and of \ff{eq:hypg}.\\
The application of \cite[Theorem 2.3]{fw16} immediately yields, in particular\footnote{We remark that \cite[Theorem 2.3]{fw16} establishes the existence of a normal form, which implies, for the class of systems considered, the perpetual stability as a consequence.}, the following
\begin{cor} Every solution of \ff{eq:hamiltonian} with initial condition $I(0) \in G$ is perpetually stable. 
\end{cor}
The paper is organized as follows: in \cref{sec:three} we cast $\psi^{(0)}$ into a suitable (integrable) form, i.e. $H_0$, in a neighbourhood of the elliptic point by means of a real-analytic transformation of variables $\mathcal{T}$. This is done by using standard methods, the integrability and the time-independence of the ``unperturbed'' travelling wave model, \cite{kw89}. Subsequently, (\cref{sec:five}) it is shown how the disturbances given by a non-flat bottom can be treated in a perturbative setting. More precisely, the implicit equation that would be obtained from \ff{eq:sw} by writing $\psi=\psi^{(0)}+\tilde{\psi}$, with $\tilde{\psi}=\mathcal{O}(\mu)$, is rewritten in terms of an infinite hierarchy of explicit equations by means of a classical formal expansion, see \ff{eq:expansion} and by using the \emph{ansatz} \ff{eq:class} for the  elements of the expansion. This is a well known approach, see e.g. \cite{gen04}. Each equation of the obtained hierarchy is shown to be resolvable in a constructive way. \\
In the final part of \cref{sec:five}, the convergence of the above described perturbative scheme is discussed in a suitable real-analyticity class. In this way, the constructed solution is mapped, through $\mathcal{T}$, to a function that will be the play the role of perturbation, i.e. $H_1$, and this will complete the construction of the Hamiltonian \ff{eq:hamiltonian}.\\
Section \ref{sec:five} uses some technical tools discussed separately in \cref{sec:four} for reader's convenience.


\section{Normal form of the unperturbed problem}\label{sec:three}
\begin{lem}\label{lem:norfor} There exists a set of action-angle variables $(I,\ph)$ and a canonical transformations $\ml{T}$ casting the streamfunction $\psi^{(0)}$, in a neighbourhood of the elliptic point $(p_e^+,q_e)$, into the form $H_0=H_0(I)$ 
\end{lem}
\proof
Let us recall \ff{eq:equilibria}. In order to study the system in a neighbourhood of the elliptic equilibrium $(p_e^+,q_e)$, let us consider the composition of the following canonical maps 
\[
\begin{array}{rcl}
\de \mathcal{C}': &  \quad & (p,q)=\left(p'+p_e^+ , q'\right)\\
\de \mathcal{C}'': &  \quad & (p',q')=\left(\ka p'', \ka^{-1} q''\right)\\
\end{array}
\]
casting the Hamiltonian (\ref{eq:ht}) in the form
\beq{eq:hs}
H_0'':= \tilde{H}_0 \circ \mathcal{C}' \circ \mathcal{C}''   = \sigma p'' - \frac{\sigma}{\ka \tilde{m}} \sin ( \ka p'') \cos(q'') - \frac{\sigma}{\ka \tilde{m}} \cos(\ka p'') \cos(q'')  \sqrt{\left(\frac{A \ka \tilde{m}}{\sigma} \right)^2-1}  \mx{.}
\eeq
Then we finally consider the following transfomation
\[  \mathcal{C}''':  \quad  (p'',q'') = ((\ka \tilde{m}) ^{-1} P,Q)  \mx{.}\]
The equations of motion in the new variables are still canonical. In fact, by setting $\lambda:=\sqrt{ \sigma^{-2} (A \ka \tilde{m})^2-1} \in \RR^+$ by \ff{eq:condparam}, their Hamiltonian reads as 
\beq{eq:finham}
H_0'''(P,Q)=\sigma [P- (\sin P + \lambda \cos P) \cos Q ] \mx{.}
\eeq
Note that the elliptic equilibrium is now $(P^*,Q^*)=(0,0)$. A Taylor expansion of (\ref{eq:finham}) centred at the origin yields
\beq{eq:hamnor} 
H_0'''(P,Q) = \frac{\sigma \lambda}{2} (P^2+Q^2)+ \mathcal{O}_3(P,Q) \mx{.} 
\eeq
It is well known, see e.g. \cite{arn}, that there exists a canonical transformation $\mathcal{B}$ casting the Hamiltonian (\ref{eq:hamnor}) into the following \emph{Birkhoff normal form}
\beq{eq:birnorfor}
H_0^{\mathcal{B}}(P,Q)=\frac{\sigma \lambda}{2}(P^2+Q^2)+\sum_{k \geq 1} \alpha_k (P^2+Q^2)^{2k} \mx{,}
\eeq 
for bounded $\{\alpha_k\} \in \RR$. The transformation is analytic in some complex neighbourhood of the origin, i.e. in a set $D_{R}:=\{(P,Q) \in \CC^2: |P|,|Q|< R\}$ with some positive $R=\mathcal{O}(1)$, see e.g. \cite{gi12}.\\
The last step consists in the use of the well known transformation 
\beq{eq:trasfp}
\mathcal{P}: \qquad (P,Q)=\left(\sqrt{\frac{2 I}{-\sigma \lambda}} \cos \ph , \sqrt{\frac{2 I}{-\sigma \lambda}} \sin \ph \right) \mx{.}
\eeq
Hence the required change of variable is  
\beq{eq:t}
\mathcal{T}:(x,y) \rw (I,\ph), \qquad \mathcal{T}:=\ml{G} \circ \ml{C}' \circ \ml{C}'' \circ \ml{C}''' \circ \ml{B} \circ \ml{P} \mx{,}  
\eeq
where $\mathcal{G}$ has been defined in \ff{eq:galileian}. Note that the interval $G$ mentioned in the statement of Theorem \ref{thm} is given by the values $\Re (P^2+Q^2)/2$ as from (\ref{eq:trasfp}) and condition $|P|,|Q|<R$. \\
Given the real-analyticity of the maps of the composition \ff{eq:t}, the map $\ml{T}$ will be real-analytic too. This implies that the same property holds for $H_0$, being the image through $\ml{T}$ of $\psi^{(0)}$, which is clearly real-analytic.


\section{Intermezzo: some technical tools}\label{sec:four}
\subsection{Estimates on a second order ODE}
\label{sec:Intermezzo}
Given a function $f:\mathcal{S}_{\hrho} \rw \CC$ we denote with $\drx{f}(x)$ and $\ddrx{f}(x)$ the derivatives $f'(x)$ and $f''(x)$, respectively.   
\begin{prop}\label{prop:ode}
Let $R:\mathcal{S}_{\hrho} \rw \CC$ such that $|R(x)| \leq M \exp(- \nu |x|)$ with $\nu, \hat{\rho} \in (0,1]$ and consider the following differential equation
\beq{eq:odemodel}
\ddrx{B} (x)-2 i \alpha \drx{B}(x) - \beta_m^2 B(x) = R(x) \mx{,}
\eeq 
where $\alpha>0$ and $\{\beta_m^2\}_{m \in \ZZ} \subset \RR^+$. Setting
\beq{eq:deltam}
\delta_m:=\sqrt{\beta_m^2-\alpha^2} \in \CC \mx{,}
\eeq
let us suppose that $\delta_m \neq 0$ at least for some $m$, then define $
\delta_-:=\inf_{m \in \ZZ}\{|\delta_m|:\beta_m^2 \neq \alpha^2\}$ and $\delta_+:=\sup_{m \in \ZZ}\{|\delta_m|:\beta_m^2 <\alpha^2\}$ or $\delta_+:=0$ if $\beta_m^2 \geq \alpha^2$ for all $m$. Then the following bounds hold
\beq{eq:solode}
|B(x)| \leq \frac{\mathfrak{G} M}{\nu^3 (1+|\delta_m|)}  e^{ -\frac{\nu |x|}{2}}, 
\qquad |\drx{B}(x)| \leq \mathfrak{G} M  e^{-\frac{\nu |x|}{2}}, \qquad \forall x \in \mathcal{S}_{\hrho} \mx{,}
\eeq
where $\mathfrak{G}:=32 (2+\alpha) e^{(\alpha+2\nu+|\delta_+|)}  \delta_-^{-1} \nu^{-3}$.
\end{prop}
\begin{rem}
The key feature of this statement consists in the possibility to provide an explicit estimate for $|\drx{B}(x)|$. This avoids the use of a Cauchy bound and the consequent domain restriction, playing a key role in the convergence argument of Section \ref{sec:five}. \\
The possibility to consider a sequence $\{\beta_m\}$ instead of $\beta \in \RR$ enable us to provide the uniform estimates \ff{eq:solode}, being $\delta_-$ a lower bound for $\delta_m>0$.
\end{rem}
The proof uses some ideas from \cite[Proposition 2]{fw14}. 
\proof We shall replace $\beta_m, \delta_m$ with $\beta, \delta$ for simplicity throughout the proof. The cases $\beta>\alpha$, $\beta=\alpha$ and $\beta < \alpha$ will be examined separately. The general solution of equation (\ref{eq:odemodel}) is
\beq{eq:solutionode}
B(x)=
\begin{cases}
\de e^{i \alpha x}(K_1+x K_2) + \int_0^x (x-y)e^{i \alpha (x-y)} R(y) dy & \mx{if} \, \beta = \alpha \\ 
\de e^{i \alpha x} (K_1 e^{\delta x} + K_2 e^{-\delta x})+\frac{1}{\delta} 
\int_0^x e^{i \alpha (x-y)} \sinh(\delta (x-y)) R(y)dy & \mx{otherwise}
\end{cases}
\eeq
\subsubsection*{Case $\beta>\alpha$.} Hence $\delta \in \RR^+$.
Denote $x=:\xi+i \eta \in \mathcal{S}_{\hrho}$. By using their path-independence, we write the integral appearing in (\ref{eq:solutionode}) as
\[
\int_0^x e^{-i \alpha y} e^{\pm \delta y} R(y) dy = \int_0^{\xi} e^{-i \alpha \xi'} 
e^{\pm \delta \xi'} R(\xi') d \xi' + i e^{-i \alpha \xi} e^{\pm \delta \xi} \int_0^{\eta} e^{\alpha \eta'} e^{\pm i \delta \eta'} R(\xi+i \eta') d \eta' \mx{.}
\] 
Now we choose 
\beq{eq:incondode}
K_1:=-\frac{1}{2 \delta} \int_{0}^{+ \infty} e^{-(i \alpha +\delta)y}  R(y) dy , \qquad K_2=0 \mx{.} 
\eeq
Note that $|K_1| < + \infty$. In this way, by setting 
\beq{eq:itilde}
\mathcal{I}_1:=e^{(i \alpha +\delta)x}\left(K_1+\frac{1}{2 \delta} \int_0^x e^{-i \alpha y} e^{-\delta y} R(y) dy \right), \qquad \mathcal{I}_2:=-\frac{e^{(i \alpha-\delta)x}}{2 \delta} \int_0^x e^{-i \alpha y} e^{\delta y} R(y) dy \mx{,}
\eeq
we have $B(x)=\mathcal{I}_1+\mathcal{I}_2$ where, in particular
\[
\mathcal{I}_1 = \frac{e^{(i \alpha+\delta)\xi}e^{(i \delta- \alpha) \eta}}{2 \delta} \left[ i e^{-(i \alpha  + \delta) \xi} \int_0^{\eta} e^{(\alpha-i \delta) \eta'}  R(\xi+i \eta') d \eta' - \int_{\xi}^{+\infty}e^{-(i \alpha +\delta)\xi'} 
 R(\xi') d \xi' \right] \mx{.}
\] 
Hence, by hypothesis
\beqar{eq:i1}
|\mathcal{I}_1|& \leq & \de \frac{M e^{\delta \xi -\alpha \eta }}{2 \delta}
\left( e^{-\delta \xi} \int_0^{\eta} e^{\alpha \eta'-\nu \xi} d \eta'+
\int_{\xi}^{+ \infty} e^{-(\delta+\nu)\xi'} d \xi'\right) \\
& \leq & \de  \frac{M e^{(\alpha + \nu) \hrho }}{2 \delta} \left( |\eta| e^{\alpha |\eta|-\nu \xi} + \frac{e^{-\nu \xi}}{\delta+\nu}\right) \\
& \leq & \de M e^{(\alpha + \nu) \hrho } (2 \delta)^{-1} \left[ e + (\delta+\nu)^{-1} \right] e^{-\nu \xi}\\
& \leq & \de 2 M e^{(\alpha + \nu) \hrho}(\nu \delta)^{-1} e^{- \nu \xi} \mx{.} 
\eeqar
On the other hand, 
\beqar{eq:i2}|\mathcal{I}_2| & \leq & \de \frac{M e^{-(\alpha \eta+\delta \xi)}}{2 \delta } \left( e^{\delta \xi} e^{\nu (\rho-\xi)} \int_0^{\eta}  d \eta' +\int_0^{\xi} e^{(\delta-\nu)\xi'}d \xi' \right) \\
& \leq & \de \frac{M e^{(\alpha+\nu)\hrho}}{2 \delta } \left( e^{1-\nu \xi} +  \frac{e^{-\delta \xi}-e^{-\nu \xi}}{\nu -\delta}  \right) \\
& \leq & \de \frac{M e^{(\alpha+\nu)\hrho}}{2 \delta } \left( e+ \xi e^{-\nu \xi/2}\right)e^{-\nu \xi/2}  \\
& \leq & 2 M e^{(\alpha+\nu)\hrho} (\nu \delta)^{-1} e^{- \nu \xi/2} \mx{.}
\eeqar
By using the elementary inequality $e^{- \nu \xi} \leq e^{-\nu |x|} e^{\nu \hrho}$ valid for all $x \in \mathcal{S}_{\hrho}$,  
\beq{eq:b}
|B(x)| \leq |\mathcal{I}_1|+ |\mathcal{I}_2| \leq  4 M e^{(\alpha+ 2 \nu) \hrho} ( \nu\delta)^{-1} e^{- \nu |x|/2} \mx{.}
\eeq
As for the first derivative of $B(x)$, we have
\[
\drx{B}(x)=i \alpha B(x)+\delta (\mathcal{I}_1-\mathcal{I}_2) + \delta^{-1} R(x) \mx{,} 
\]
hence, by (\ref{eq:b}), (\ref{eq:i1}), (\ref{eq:i2}) and by hypothesis, respectively
\beq{eq:bprime}
|\drx{B}(x)|  \leq   M (4 \alpha (\nu \delta)^{-1} e^{(\alpha+2 \nu)\hrho}+4 \nu^{-1} e^{(\alpha+\nu)\hrho} + \delta^{-1} ) e^{-\nu |x|} \leq 4 M e^{(\alpha+2 \nu)\hrho} (2+\alpha) \nu^{-1} e^{- \nu |x|/2} \mx{.}
\eeq
\subsubsection*{Case $\beta<\alpha$.}
In this case $\delta = i \gamma $ with $\gamma>0$. We shall choose $K_1$ and $\mathcal{I}_1$ as in \ff{eq:incondode} and \ff{eq:i1}, respectively, and 
\[
\mathcal{I}_2:=e^{i(\alpha-\gamma)x} \left(K_2- \frac{1}{2 \delta} \int_0^x e^{-i(\alpha-\delta)y}R(y) dy\right) , \qquad
K_2:=\frac{1}{2 \delta} \int_0^{+\infty} e^{-i(\alpha+\gamma)y}R(y)dy \mx{,}
\]
note that $|K_2|<\infty$ by hypothesis on $R(x)$. The procedure is similar to \ff{eq:i1} and \ff{eq:i2} simply replacing $\delta$ with $i \gamma$. The only difference lies in the term $\exp(\gamma \eta)$ which can be bounded by $\exp(|\delta_+|)$. This yields 
\beq{eq:i1i2betalessalpha}
|\mathcal{I}_1| \leq 2 M e^{(\alpha+\nu)\hrho} (\gamma \nu)^{-1} e^{-\nu \xi/2}, \qquad 
|\mathcal{I}_2| \leq 2 M e^{[(\alpha+\nu)\hrho+|\delta_+|]} (\gamma \nu)^{-1} e^{-\nu \xi/2}, 
\eeq
then
\beq{eq:bbetalessalpha}
|B(x)| \leq 4 M e^{(\alpha+2 \nu)\hrho}(\nu \gamma)^{-1}e^{-\nu |x|/2}, \qquad 
|\drx{B}(x)| \leq 4 M e^{[(\alpha+\nu)\hrho+|\delta_+|]} \nu^{-1} e^{-\nu |x|/2} \mx{.}
\eeq
\subsubsection*{Case $\beta=\alpha$.} This case ($\delta=0$) can be treated in a similar way. We only mention that one can choose
\[
K_1=\int_0^{+\infty} \xi' e^{-i \alpha \xi'}R(\xi') d \xi', \qquad 
K_2=-\int_0^{+\infty} e^{-i \alpha \xi'} R(\xi') d \xi' \mx{,}
\]
obtaining
\beq{eq:estimatexaub}
|B(x)|, |\drx{B}(x)| \leq 16 M (1+\alpha) \nu^{-3} e^{(\alpha+2 \nu)\hrho} e^{-\nu |x|/2} \mx{.}
\eeq
The required bounds follow immediately by collecting those obtained in \ff{eq:b}, \ff{eq:bprime}, \ff{eq:bbetalessalpha} and \ff{eq:estimatexaub} and using the fact that $|\delta_m|^{-1} \leq 2[\delta_-(1+|\delta_m|)]^{-1}$. 
\endproof

\subsection{A bound on certain brackets of two sequences}

Let us now consider two sequences $\mathbf{f}:=\{f_n(x)\}_{n \in \ZZ }$ and $\mathbf{g}:=\{g_n(x)\}_{n \in \ZZ }$, then define (formally),  for all $m \in \ZZ$, 
\beq{eq:brack}
[\mathbf{f},\mathbf{g}]_m:=\sum_{l \in \ZZ} l(f_l \drx{g}_{m-l}-g_l \drx{f}_{m-l}) \mx{,}
\eeq
which is a sequence of functions indexed by $m$.

\begin{prop}\label{prop:prod} Let $\{f_l(x)\}$ and $\{g_l(x)\}$ be two sequences, with $f_l(x),g_l(x):\ml{S}_{\hat{\rho}} \rw \CC$ for all $l \in \ZZ$, such that, for all $\hat{\rho}<\rho$ and all $l \in \ZZ \setminus\{0\}$,  
\beq{eq:prophyp}
|f_l(x)|\leq \frac{K(x)}{|l|} e^{-|l|\hat{\rho}}, \qquad |g_l(x)| \leq L(x) e^{-|l|\rho} \mx{.}
\eeq
Then, denoted $\delta:=\rho-\hat{\rho}>0$, one has
\beq{eq:bound}
|[\mathbf{f},\mathbf{g}]_m| \leq 8 K(x) L (x) \delta^{-3} e^{-|m|\hat{\rho}} \mx{.}
\eeq
\end{prop}
\proof Along the lines of \cite[Lemma 4.1]{gi03}, by a Cauchy estimate
\[
|\drx{g}_l(x)| \leq \delta^{-1} |g_l(x)|, \quad \forall x \in \mathcal{S}_{\hat{\rho}} \mx{,}
\] 
hence
$$
\begin{array}{rcl}
|[\mathbf{f},\mathbf{g}]_m| & \leq & \de K(x) L(x) \delta^{-1} \sum_{l \neq 0}  \left[ e^{-|l|\hat{\rho}}e^{-|m-l|\rho}+|l|e^{-|l|\rho e^{-|m-l|\hat{\rho}}}\right]\\
& = & \de K(x) L(x) \delta^{-1} \sum_{l \neq 0} \left[ e^{-|m-l|\delta}+|l|e^{-|l|\delta} \right] e^{-(|l|+|m-l|)\hat{\rho}} \\
& \leq & \de 4 K(x) L(x) \delta^{-1} \left[\delta^{-1}+\delta^{-2}\right] e^{-|m|\hat{\rho}} \mx{,}
\end{array}
$$  
where we have used the elementary inequality $\sum_{l \geq 1} l^a \exp(-l \delta) \leq 2 \delta^{-(1+a)}$, for all $\delta \in (0,1/2]$ and $a=0,1$.
\endproof
\begin{rem} We anticipate here that the presence of the term $|l|^{-1}$ in the first of (\ref{eq:bound}) will play a crucial role in the proof of Lemma \ref{lem:convergence}. Without its contribution it would have not been possible to ``preserve'' the term $\exp(-|m| \hat{\rho})$ in (\ref{eq:bound}). As it is well known in the real-analytic functions context, see e.g. \cite{gi03}, the absence of this term would have implied a domain restriction, namely, the $\exp(-|m| \hat{\rho})$ above, would have been replaced by a term of the form $C d^{-1} \exp(-(1-d)|m|\hat{\rho})$, for some $d>0$ (restriction) and some constant $C>0$. \\ 
As a standard procedure in perturbation theory one chooses a sequence of $d$'s  e.g. $d_j \sim j^{-2}$ in such a way the final domain  $\Pi_{j \geq 1} (1-d_j)\rho$ is non-trivial. Unfortunately, the accumulation of these $d_j^{-1}$ would invalidate the convergence argument of Lemma \ref{lem:convergence} as can be easily seen from \ff{eq:ep} and \ff{eq:last}. This kind of obstruction to Cauchy's majorants method is a well known phenomenon and it cannot be overcome unless the above mentioned ``artificial'' \emph{small divisors} are controlled with some specialized argument, see for instance \cite{gi12}.\\
\end{rem}



\section{The perturbative setting}\label{sec:five}
\subsection{Formal scheme}
Let us consider equation (\ref{eq:sw}) in which the bathymetry has been chosen of the form (\ref{eq:bath}) and write 
\beq{eq:psisum}
\psi(x,y,t)=\psi^{(0)}(x,y,t)+\tilde{\psi}(x,y,t) \mx{,}
\eeq
then expanding\footnote{A similar expansion of the streamfunction, despite with respect to a different parameter, was already considered in \cite{rhi73}.}, formally 
\beq{eq:expansion}
\tilde{\psi}(x,y,t)=\sum_{j=1}^{\infty} \lambda^j \psi^{(j)}(x,y,t)  \mx{,}
\eeq
where $\lambda$ is an auxiliary parameter which will be thought to be equal to one (\emph{book-keeping parameter}, see e.g. \cite{efthy}). On the oher hand, let $g:=\mu \tilde{g}$ we get
\[
J(\psi,\lambda g-d)=\sum_{j=1}^{\infty }\lambda^j \left[\pl_x \psi^{(j-1)} \pl_y g-\pl_y \psi^{(j-1)} \pl_x g \right]= \sum_{j=1}^{\infty }\lambda^j J(g,\psi^{(j-1)}) \mx{.}
\]
In this setting, equation (\ref{eq:sw}) is equivalent to the following hierarchy of recursive equations (obtained by balancing the powers of $\lambda$) 
\beq{eq:hier}
\pl_t (\Delta \psi^{(j)}-F \psi^{(j)})+\mathcal{F} \pl_x \psi^{(j)}= J(g,\psi^{(j-1)})
\mx{,}
\eeq
for all $j \geq 1$, where the order zero is given by (\ref{eq:swzero}).\\
The following statement shows the possibility to solve, at least formally, the  equations of the hierarchy up to an arbitrarily high order in a constructive way 
\begin{prop}\label{prop:formal} For all $j \geq 1$, it is possible to determine a sequence of functions $\{b_{m,n}^{(j)}(x)\}$ satisfying
\begin{subequations} 
\begin{align}
b_{-m,n}^{(j)}(x) & =-b_{m,n}^{(j)}(x) \label{eq:simone}\\
b_{\mp m,n}^{(j)}(x) & = \bar{b}_{ \pm m, -n}^{(j)}(x) \label{eq:simtwo} 
\end{align} 
\end{subequations}
in particular $b_{0,n}^{(j)}(x) \equiv 0$ ($\bar{z}$ denotes the complex-conjugate of $z \in \CC$) and such that the elements
\beq{eq:class}
\psi^{(j)}(x,y,t)=\sum_{(m,n) \in \mathcal{J}} b_{m,n}^{(j)}(x) e^{i (m y+\sigma(n)t)} \mx{,}
\eeq
where $\mathcal{J}:=\{(m,n):m \in \ZZ \setminus \{0\}, \, n=\pm \kappa \}$, are real solutions of the hierarchy (\ref{eq:hier}) for all $(x,y) \in \mathcal{C}$ and all $t \in \RR$, under the boundary conditions (\ref{eq:boundary}).
\end{prop} 
\proof Let us preliminarily observe that, from the standard theory of Fourier series, the conditions (\ref{eq:simone}) and (\ref{eq:simtwo}) hold if, and only if, $\psi^{(j)}(x,y,t)$ is odd in $y$ and real, respectively. In particular, (\ref{eq:simone}) implies that the conditions (\ref{eq:boundary}) are satisfied.\\
Note that  (\ref{eq:psizero}) can be written in the form (\ref{eq:class}) with 
\beq{eq:bzero}
b_{m,n}^{(0)}(x)= 
\begin{cases}
A (4 i)^{-1} \sign(m) \exp(i n x) & \mx{if} \quad (m,n)=(\pm \tilde{m},\pm \ka) \\
0 & \mx{otherwise}
\end{cases} \mx{.}
\eeq
Hence (\ref{eq:simone}) holds, being (\ref{eq:simtwo}) obvious.\\
Then we can suppose that (\ref{eq:class}) satisfies the equations (\ref{eq:hier}) with $b_{m,n}^{(j)}(x)$ recursively determined and satisfying (\ref{eq:simone}) and (\ref{eq:simtwo}) up to a level $j-1$ and proceed by induction. Let us write 
\[
g(x,y)=\sum_{l \in \ZZ \setminus \{0\}} g_l(x) e^{i l y} \mx{,}
\]
where, by hypothesis, $g_0(x) \equiv 0$ and  
\beq{eq:propg}
g_{-l}(x)=-g_{l}(x), \quad g_{-l}(x)=\bar{g}_{l}(x) \mx{,} 
\eeq
Hence, 
\[
\begin{array}{rcl}
J(g,\psi^{(j-1)}) & = & \de i \sum_{n = \pm \kappa} e^{i \sigma(n) t } 
\left[ \left( \sum_l \drx{g}_l e^{i l y}\right)  \left( \sum_m m b_{m,n}^{(j-1)} e^{i m y} \right)-
 \left( \sum_l l g_l e^{i l y}\right)\left( \sum_m \drx{b}_{m,n}^{(j-1)} e^{i m y} \right)
\right]\\
& = & \de i \sum_{n = \pm \kappa} e^{i \sigma(n) t } \sum_{m \neq 0} \left[ \sum_{k'} k' \left( b_{k',n}^{(j-1)} \drx{g}_{m-k'}- g_{k'} \drx{b}_{m-k',n}^{(j-1)} \right)\right]\\
& = & i  \de \sum_{(m,n) \in \mathcal{J}} [\mathbf{b}_{\cdot,n}^{(j-1)},\mathbf{g}]_m e^{i (m y+\sigma(n) t) } \mx{,} 
\end{array}
\]
where the formula of the product \emph{\`{a} la Cauchy}, the definition (\ref{eq:brack}) and the fact that $ [\mathbf{b}_{\cdot,n}^{(j-1)},\mathbf{g}]_0=0$ due to (\ref{eq:simone}) and (\ref{eq:propg}) have been used.\\ 
On the other hand, by substituting $\psi^{(j)}$ as in (\ref{eq:class}) in the l.h.s. of (\ref{eq:hier}), we get that the latter is satisfied if the following linear differential equation 
\beq{eq:ode}
\sigma(n) \ddrx{b}_{m,n}^{(j)} - i \mathcal{F} \drx{b}_{m,n}^{(j)} -\sigma(n) (F + m^2) b_{m,n}^{(j)} = [\mathbf{b}_{\cdot,n}^{(j-1)},\mathbf{g}]_m \mx{,}
\eeq
holds true for all $(m,n) \in \mathcal{J}$. Note that $[\mathbf{b}_{\cdot,n}^{(j-1)},\mathbf{g}]_m $ is a function of $x$, known by hypothesis. By dividing both sides of (\ref{eq:ode}) by $\sigma(n) < 0 $, we have that equation (\ref{eq:ode}) is of the form (\ref{eq:odemodel}). Hence $b_{m,n}^{(j)}(x)$ can be uniquely determined by (\ref{eq:solutionode}). In addition, $b_{m,n}^{(j)}(x)$ satisfies the relations (\ref{eq:simone}) and (\ref{eq:simtwo}) if $[\mathbf{b}_{\cdot,n}^{(j-1)},\mathbf{g}]_m$ does. The latter property is a straightforward check under the hypothesis (\ref{eq:propg}).\\  
This completes the proof and the formal resolvability of the hierarchy (\ref{eq:hier}).
\endproof
\subsection{Convergence}
In this section, the convergence of the formal scheme built in the previous section is addressed, and it can be stated as in the following 
\begin{lem}\label{lem:convergence} Set $\beta_m^2:=F+m^2$ and let $m^*$ be such that $\delta_{m^*}=\delta_-$ as defined in Proposition \ref{prop:ode}. \\  
If $\mu$ satisfies
\beq{eq:mustar}
\mathcal{L}(\mu):=\frac{32 \mu \ml{M} \mathfrak{G}(|m^*|+1) }{\sigma \rho^{2}} \leq \frac{1}{2} \mx{,}  
\eeq
then $\tilde{\psi}$, as defined in \ff{eq:expansion}, is a real-analytic function on $\mathcal{C}_{\rho/4}$, satisfying
\beq{eq:finalestimate}
|| \tilde{\psi}(x,y,t)||_{\rho/4} \leq C(\mu) e^{-\nu  \frac{|x|}{2}} \mx{.} 
\eeq
where $C(\mu)$ is a $O(\mu)$ constant.
\end{lem}
Before proceeding with the proof, let us observe that the statement above implies Theorem \ref{thm}. The term $H_0(I)$ is given by Lemma \ref{lem:norfor}. On the other hand, by setting  
\[
H_1(x,y,t):=\tilde{\psi} \circ \mathcal{T} \mx{,}
\]
where $\mathcal{T}$ has been defined in \ff{eq:t} we clearly have $H_1$ is $\mathcal{O}(\mu)$ by \ff{eq:finalestimate}. Furthermore, as $\tilde{\psi}$ is real-analytic on  $\mathcal{C}_{\rho/4}$ and $\mathcal{T}$ is a real-analytic map then $H_1$ will be real-analytic as well, more precisely in $\mathcal{U}_{\tilde{\gamma}}$ for some $\tilde{\gamma}>0$. Finally, it is clear that the transformation $\mathcal{G}$ itself, produces a linear growth in time of $\Re x$. Hence, the function $\tilde{\psi}$ which varies with a general dependence in $x$, see \ff{eq:class}, will be (at least in general) \emph{aperiodic} in time and, by \ff{eq:finalestimate} and \ff{eq:galileian}, it will decay in time as in (\ref{eq:decayintime}). \endproof
\proof
Following the classical approach of the Cauchy majorants, the strategy consists in showing, by induction, that there exists a suitable infinitesimal sequence $\{\epsilon_j\}_{j \in \NN}$ such that  
\beq{eq:estimatebnm}
|m||b_{m,n}^{(j)}(x)|, |\drx{b}_{m,n}^{(j)}(x)|  \leq 
\begin{cases}
\epsilon_0 e^{-|m|\frac{\rho}{2}} & \mx{if} \, j=0\\ 
\epsilon_j e^{-|m|\frac{\rho}{2}}e^{-\nu \frac{|x|}{2}} & \mx{otherwise}  
\end{cases} \mx{.}
\eeq
First of all we note that, only a finite number of $b_{m,n}^{(0)}$ are different from zero, i.e., $b_{m,n}^{(0)}(x) \equiv 0 $ for all $|m| \geq \tilde{M}$. Furthermore, by (\ref{eq:bzero})
\[
\tilde{\epsilon}:=\max_{\substack{|m| \leq \tilde{M} \\ n = \pm \kappa }}\{ \sup_{x \in \mathcal{S}_{\rho}}|b_{m,n}^{(0)}(x)|,  \sup_{x \in \mathcal{S}_{\rho}}|\drx{b}_{m,n}^{(0)}(x)|\} = \frac{\kappa A}{4} \mx{.} 
\]
By setting $\epsilon_0:=\tilde{M} \tilde{\epsilon}$ the induction basis ($j=0$) is clearly true. \\
Let us suppose the statement to be true for $j$ and proceed by induction. By \ff{eq:hypg} we have 
\beq{eq:g}
|g_l(x)| \leq \epsilon e^{-|l|\rho} e^{-\nu |x|},\qquad \forall l \in \ZZ \setminus\{0\} \mx{,}
\eeq 
hence we can set $K(x):=\epsilon_j$ (we disregard here the term $\exp(-\nu|x|/2)$) and $L(x):=\epsilon \exp(-\nu |x|)$ and we use Proposition \ref{prop:prod} with $\hat{\rho}=\rho/2$. The latter yields, in particular
\[
|[\mathbf{b}_{\cdot,n}^{(j)},\mathbf{g}]_m| \leq 32 \epsilon_j \epsilon \rho^{-1} e^{-|m| \frac{\rho}{2}} e^{-\nu |x|} \mx{.}
\]
We are now able to give estimate on the solutions of equation (\ref{eq:ode}). In fact, by setting $\alpha \leftarrow \mathcal{F}/(2 \sigma)$ and $\beta_m$ as in the statement, the equation (\ref{eq:ode}) is of the form (\ref{eq:odemodel}) with $R:=\sigma^{-1}[\mathbf{b}_{\cdot,n}^{(j)},\mathbf{g}]_m $. \\
Hence, by Proposition \ref{prop:ode} (with $\hat{\rho} \leftarrow \rho/2$) we get 
\[
|m| |b_{m,n}^{(j+1)}(x)|, |\drx{b}_{m,n}^{(j+1)}(x)| \leq 32 \mathfrak{G} (|m^*|+1) \sigma^{-1} \rho^{-2} \epsilon \epsilon_j e^{-|m| \frac{\rho}{2}} e^{-\nu |x|/2} \mx{,}
\]
where we have used the bound $(1+|\delta_m|)^{-1} \leq (|m^*|+1)/|m|$. \\ 
Hence recalling that $\epsilon=\mu \ml{M}$ and by setting
\beq{eq:ep}
\epsilon_{j+1}:=\mathcal{L}(\mu) \epsilon_j \mx{,}
\eeq
the inductive step is complete.\\
Furthermore, (\ref{eq:ep}) implies that $\epsilon_j:=\mathcal{L}^j(\mu) \epsilon_0$, then by (\ref{eq:estimatebnm}) and under the condition (\ref{eq:mustar}), 
\beqar{eq:last}
\norm{\tilde{\psi}}{\rho/4} & \leq & \de \sum_{j \geq 1} \sum_{(m,n) \in \mathcal{J}} |b_{m,n}\jj (x)| e^{|m|\frac{\rho}{4}} \\
& \leq & \de 2  e^{-\nu |x|/2} \sum_{j \geq 1} \epsilon_j \sum_{m \in \ZZ \setminus \{0\}} e^{-|m|\frac{\rho}{4}} \\
& \leq & \de 4 \epsilon_0 \mathcal{L}(\mu) e^{-\frac{\rho}{4}} \left[(1-e^{-\frac{\rho}{4}})(1-\mathcal{L}(\mu)) \right]^{-1}  e^{-\nu |x|/2} \mx{,}
\eeqar
which implies (\ref{eq:finalestimate}). 
\endproof

\subsection*{Acknowledgements}
This research was supported by ONR Grant No. N00014-01-1-0769.\\
The author thanks Prof. S. Wiggins for stimulating discussions. The author is also grateful to Prof. L. Biasco for very useful comments on a preliminary version of this manuscript.

\bibliographystyle{alpha}
\bibliography{Geo.bib}

\begin{thebibliography}{IAKN14}

\bibitem[Ber07]{be03}
M.~Berti.
\newblock {\em Nonlinear Oscillations of Hamiltonian PDEs}.
\newblock Progress in Nonlinear Differential Equations and Their Applications.
  Birkh{\"a}user Boston, 2007.

\bibitem[Cha07]{ch07}
F.~Chazel.
\newblock Influence of bottom topography on long water waves.
\newblock {\em M2AN Math. Model. Numer. Anal.}, 41(4):771--799, 2007.

\bibitem[CLS12]{cr12}
W.~Craig, D.~Lannes, and C.~Sulem.
\newblock Water waves over a rough bottom in the shallow water regime.
\newblock {\em Ann. Inst. H. Poincar\'e Anal. Non Lin\'eaire}, 29(2):233--259,
  2012.

\bibitem[CS10]{cr10}
W.~Craig and C.~Sulem.
\newblock Asymptotics of surface waves over random bathymetry.
\newblock {\em Quart. Appl. Math.}, 68(1):91--112, 2010.

\bibitem[Eft12]{efthy}
C.~Efthymiopoulos.
\newblock Canonical perturbation theory, stability and diffusion in
  {H}amiltonian systems. applications in dynamical astronomy.
\newblock In C.~Giordano P.~Cincotta and C.~Efthymiopoulos, editors, {\em
  Proceedings of the 3rd La Plata School on Astronomy and Geophysics,
  Association of Astronomy of Argentina}, pages 1--144., 2012.

\bibitem[FW14]{fw14}
A.~Fortunati and S.~Wiggins.
\newblock Persistence of diophantine flows for quadratic nearly integrable
  hamiltonians under slowly decaying aperiodic time dependence.
\newblock {\em Regular and Chaotic Dynamics}, 19(5):586--600, 2014.

\bibitem[FW16]{fw16}
A.~Fortunati and S.~Wiggins.
\newblock Negligibility of small divisor effects in the normal form theory for
  nearly-integrable hamiltonians with decaying non-autonomous perturbations.
\newblock {\em Celestial Mechanics and Dynamical Astronomy}, 125(2):247--262,
  2016.

\bibitem[Gio03]{gi03}
A.~Giorgilli.
\newblock Exponential stability of {H}amiltonian systems.
\newblock In {\em Dynamical systems. {P}art {I}}, Pubbl. Cent. Ric. Mat. Ennio
  Giorgi, pages 87--198. Scuola Norm. Sup., Pisa, 2003.

\bibitem[Gio12]{gi12}
A.~Giorgilli.
\newblock On a {T}heorem of {L}yapounov.
\newblock {\em Rendiconti dell'Istituto Lombardo Accademia di Scienze e
  Lettere, Classe di Scienze Matematiche e Naturali}, 146:133--160, 2012.

\bibitem[GM04]{gen04}
G.~Gentile and V.~Mastropietro.
\newblock Convergence of lindstedt series for the nonlinear wave equation.
\newblock {\em Communications on Pure and Applied Analysis.}, 3:509--514, 2004.

\bibitem[IAKN14]{arn}
A.~Iacob, V.I. Arnol'd, V.V. Kozlov, and A.I. Neishtadt.
\newblock {\em Dynamical Systems III}.
\newblock Encyclopaedia of Mathematical Sciences. Springer Berlin Heidelberg,
  2014.

\bibitem[KW87]{kn87}
E.~Knobloch and J.~B. Weiss.
\newblock Chaotic advection by modulated traveling waves.
\newblock {\em Phys. Rev. A}, 36:1522--1524, Aug 1987.

\bibitem[MW98]{ma98}
N.~Malhotra and S.~Wiggins.
\newblock Geometric structures, lobe dynamics, and {L}agrangian transport in
  flows with aperiodic time-dependence, with applications to {R}ossby wave
  flow.
\newblock {\em J. Nonlinear Sci.}, 8(4):401--456, 1998.

\bibitem[Ped12]{psky}
J.~Pedlosky.
\newblock {\em Geophysical Fluid Dynamics}.
\newblock Springer Study Edition. Springer New York, 2012.

\bibitem[Pie91]{pi91}
R.~T. Pierrehumbert.
\newblock Chaotic mixing of tracer and vorticity by modulated travelling rossby
  waves.
\newblock {\em Astrophys. Fluid Dyn}, 58:285--319, 1991.

\bibitem[RB73]{rhi73}
P.~Rhines and F.~Bretherton.
\newblock Topographic rossby waves in a rough-bottomed ocean.
\newblock {\em Journal of Fluid Mechanics}, 61(3):583–607, 1973.

\bibitem[RP83]{ro83}
R.~R. Rosales and G.~C. Papanicolaou.
\newblock Gravity waves in a channel with a rough bottom.
\newblock {\em Stud. Appl. Math.}, 68(2):89--102, 1983.

\bibitem[SW06]{sam06}
R.M. Samelson and S.~Wiggins.
\newblock {\em Lagrangian Transport in Geophysical Jets and Waves: The
  Dynamical Systems Approach}.
\newblock Interdisciplinary Applied Mathematics. Springer New York, 2006.

\bibitem[Van03]{van03}
J.~Vanneste.
\newblock Nonlinear dynamics over rough topography: homogeneous and stratified
  quasi-geostrophic theory.
\newblock {\em Journal of Fluid Mechanics}, 474:299--318, 2003.

\bibitem[WK89]{kw89}
J.B. Weiss and E.~Knobloch.
\newblock Mass transport and mixing by modulated traveling waves.
\newblock {\em Phys. Rev. A}, 40:2579--2589, Sep 1989.

\end{thebibliography}

\end{document}